\documentclass[12pt]{article}

\usepackage{amsmath}
\usepackage{amsthm}
\usepackage{amssymb}
\usepackage{amsfonts}
\usepackage{bbm}
\usepackage{dsfont}
\usepackage{graphicx}
\usepackage{hyperref}

\newtheorem{theo}{Theorem}
\newenvironment{theorem}{\vspace{4mm}\begin{theo}}{\end{theo}}
\newtheorem{lem}[theo]{Lemma}
\newenvironment{lemma}{\vspace{4mm}\begin{lem}}{\end{lem}}
\newtheorem{coro}[theo]{Corollary}
\newenvironment{corollary}{\vspace{4mm}\begin{coro}}{\end{coro}}
\newtheorem{rem}[theo]{Remark}
\newenvironment{remark}{\vspace{4mm}\begin{rem}\rm}{\end{rem}}
\newtheorem{conjec}{Conjecture}

\newtheorem{prop}[theo]{Proposition}

\newtheorem{defi}{Definition}

\newtheoremstyle{citing}{}{}{\itshape}{}{\bfseries}{.}%
 { }{\thmnote{#3}}
\theoremstyle{citing}

\newcommand{\Z}{\mathbb{Z}}
\newcommand{\N}{\mathbb{N}}
\newcommand{\R}{\mathbb{R}}

\newcommand{\ind}{\scalebox{1.2}{\raisebox{-0.2mm}{$\mathds{1}$}}}
\newcommand{\abs}[1]{\left\vert #1 \right\vert}	
\newcommand{\norm}[1]{\left\Vert #1 \right\Vert}	

\renewcommand{\l}{\langle}
\renewcommand{\r}{\rangle}

\newcommand{\g}{\gamma}
\newcommand{\G}{\Gamma}

\renewcommand{\O}{\Omega}

\begin{document}

\title{Swendsen-Wang is faster than single-bond dynamics}
\author{Mario Ullrich\footnote{The author was supported 
																by the DFG GRK 1523.}\\
	Mathematisches Institut, Universit\"at Jena\\   
	email: ullrich.mario@gmail.com}



\maketitle

\begin{abstract}
We prove that the spectral gap of the Swendsen-Wang dynamics 
for the random-cluster model is larger than the spectral gap of a 
single-bond dynamics, that updates only a single edge per step. 
For this we give a representation of the algorithms on the joint 
(Potts/random-cluster) model. \\
Furthermore we obtain upper and lower bounds on the mixing time of the 
single-bond dynamics on the discrete $d$-dimensional torus of side length $L$ 
at the Potts transition temperature for $q$ large enough 
that are exponential in $L^{d-1}$, complementing a result of 
Borgs, Chayes and Tetali~[\emph{Probab. Theory Related Fields}, 152 (2012), pp. 509--557].
\end{abstract}

\section{Introduction}

This work was motivated by the recent article of Borgs, Chayes and Tetali 
``Tight bounds for mixing of the Swendsen-Wang algorithm at the Potts
transition point'' \cite{BCT}, where the authors prove upper and lower bounds 
for the 
mixing time of the Swendsen-Wang and heat-bath dynamics for the $q$-state 
Potts model on rectangular subsets of the lattice 
$\Z^d$ with periodic boundary conditions at the Potts transition temperature 
if $q$ and the side length $L$ are large enough. 
Both, upper and lower bounds, are exponential in $L^{d-1}$. 
(Their upper bounds are valid for all $q$, $L$.)
Since one can sample from the Potts model if one can do so for the 
random-cluster model, see e.g. \cite{G1}, 
we wonder if the same upper and lower bounds are valid for 
the heat-bath 
dynamics for the random-cluster model, i.e. the 
\emph{single-bond dynamics}.

In this article we give a positive answer to this question. 
For the lower bound we prove that if $G$ is an arbitrary graph, 
$p\in(0,1)$ and $q\in\N$, then 
\[
\lambda(P_{\rm SW}) \;\ge\; \lambda(P_{\rm SB}),
\]
where $\lambda(\cdot)$ denotes the spectral gap and $P_{\rm SW}$ 
(resp.~$P_{\rm SB}$) denotes the transition matrix of the 
Swendsen-Wang (resp. single-bond) dynamics as defined in 
Section~\ref{sec:alg} (see Theorem~\ref{th:main}). 
For this we represent both Markov chains by transition matrices on 
the joint (Potts/random-cluster) model. 
Estimates of the norm of the corresponding Markov operators lead to 
the result.

By this inequality and \cite{BCT} we obtain an exponential 
(in $L^{d-1}$) lower bound for the mixing time of the single-bond dynamics 
on boxes of side length $L$ in $\Z^d$ with periodic boundary condition at 
the transition point for $q$ and $L$ large enough, 
like for the Swendsen-Wang dynamics. 
The proof of the upper bound uses a result of 
Ge and {\v{S}}tefankovi{\v{c}} \cite{GeS} that provides us with 
an upper bound on the mixing time of the single-bond dynamics in 
terms of the \emph{linear-width} of a graph 
(see Section \ref{sec:torus}).

Furthermore we obtain some rapid mixing results for the Swendsen-Wang 
dynamics.
First we easily obtain rapid mixing of the 
Swendsen-Wang dynamics on trees (which was proven by several authors before, 
see e.g. \cite{CF}, \cite{BCT}), because if $T=(V,E)$ is a tree, 
the random-cluster measure is a product measure and so the single-bond dynamics 
has spectral gap $\O(1/\abs{E})$ 
(see e.g.~\cite[Lemma~12.11]{LPW}).
Hence, $\lambda(P_{\rm SW})^{-1}=\mathcal{O}(\abs{E})$ for every tree $T$, 
$p\in(0,1)$ and $q\in\N$ (see Corollary \ref{coro:tree}).\\
Additionally we get (again by \cite{GeS}) that Swendsen-Wang is rapidly mixing 
for graphs with bounded linear-width (Corollary \ref{coro:width}).

\section{The models} \label{sec:models}

Fix some $p\in(0,1)$, a natural number $q\ge1$ and a graph $G=(V,E)$ 
with finite vertex set 
$V$ and edge set $E$.

The \emph{random-cluster} (RC) \emph{model} 
(also known as the FK-model), see Fortuin and Kasteleyn \cite{FK}, 
is defined on the graph $G=(V,E)$ by its state space 
$\O_{\rm RC}=\{A: A\subseteq E\}$ and the random-cluster measure
\[
\mu(A) \;:=\; \mu^G_{p,q}(A) \;=\; 
\frac1{Z(G,p,q)}\,
\left(\frac{p}{1-p}\right)^{\abs{A}}\,q^{c(A)},
\]
where $c(A)$ is the number of connected components 
in the graph $(V,A)$, counting isolated vertices as a component, 
and $Z$ is the normalization constant that makes $\mu$ a 
probability measure.  
For a detailed introduction and related topics see \cite{G1}.
Although $\mu$ is well-defined for all $q>0$, we are only interested 
in integer values. In this case there is a tight connection to 
a model on (not necessarily proper) colorings of the vertices of the graph $G$.
The \emph{$q$-state Potts model} on $G$ at inverse temperature 
$\beta\ge0$ is defined as the set of 
possible \emph{configurations} 
$\O_{\rm P}=[q]^V$, where $[q]\,{:=}\,\{1,\dots,q\}$ 
is the set of \emph{colors} (or spins),
together with the probability measure
\vspace{1mm} 
\label{page:pi}
\[
\pi(\sigma) \;:=\; \pi^{G}_{\beta,q}(\sigma) \;=\; 
	\frac1{Z(G,1-e^{-\beta},q)}\,
	\exp\left\{\beta\,\sum_{u,v:\, u\leftrightarrow v}
	\ind\bigl(\sigma(u)=\sigma(v)\bigr)\right\}
\]
for $\sigma\in\O_{\rm P}$, where $u\leftrightarrow v$ if and only if  
$u$ and $v$ are \emph{neighbors} in $G$, i.e. $\{u,v\}\in E$, and
$Z(\cdot,\cdot,\cdot)$ is the same normalization 
constant as for the random-cluster model (see \cite[Th. 1.10]{G1}).

To describe the algorithms we also need the coupling of the 
Gibbs measure 
$\pi_{\beta,q}^G$ and the random-cluster measure $\mu_{p,q}^G$ 
of Edwards and Sokal \cite{ES}. 
Let us define 
\[
\O(A) \;:=\; \bigl\{\sigma\in\O_{\rm P}:\,\sigma(u)=\sigma(v)\;\;
	\text{ for all } \{u,v\}\in A\bigr\}, \quad A\subset E,
\]
and
\[
E(\sigma) \;:=\; \bigl\{\{u,v\}\in E:\,\sigma(u)=\sigma(v)\bigr\},
	\quad \sigma\in\O_{\rm P}.
\]
Obviously, we have for $\sigma\in\O_{\rm P}$ and $A\subset E$ that 
$\sigma\!\in\!\O(A)\Leftrightarrow A\!\subset\!E(\sigma)$. 
Let $\sigma\in\O_{\rm P}$, $A\in\O_{\rm RC}$ and $p=1-e^{-\beta}$, 
then the joint measure of $(\sigma,A)\in\O_{\rm J}:=\O_{\rm P}\times\O_{\rm RC}$ is
\[
\bar\mu(\sigma,A)\;:=\;\bar\mu_{p,q}^G(\sigma,A)
	\;=\;\frac{1}{Z(G,p,q)}\,
	\left(\frac{p}{1-p}\right)^{\abs{A}}\,\ind(A\subset E(\sigma)).
\]
The marginal distributions of $\bar\mu$ are exactly $\pi$ and 
$\mu$, respectively, and we will call $\bar\mu$ the \emph{FKES} 
(Fortuin-Kasteleyn-Edwards-Sokal) \emph{measure}.

\section{Mixing time} \label{sec:mixing}

In the following we want to estimate the efficiency of  
Markov chains. For an introduction to Markov 
chains and techniques to bound the convergence rate to the stationary 
distribution, see, e.g.,~\cite{LPW}. 
Let $P$ be the transition matrix of a Markov chain with state space 
$\Omega$ that is ergodic, i.e. irreducible and aperiodic, 
and has unique stationary measure $\pi$.
Then we define the \emph{mixing time} of the Markov chain by
\[
\tau(P) \;:=\; \min\left\{t:\; 
	\max_{x\in\O}\sum_{y\in\O}\abs{P^t(x,y)-\pi(y)}\,\le\,\frac{1}{e}\right\}.
\]
If we are considering simultaneously a \emph{family} of state spaces
$\{\O_n\}_{n\in\N}$ with a corresponding family of Markov chains 
$\{P_n\}_{n\in\N}$, we say that the chain is \emph{rapidly mixing}
for the given family if $\tau(P_n)^{-1} = \mathcal{O}(\log(|\O_n|)^C)$
for all $n\in\N$ and some $C < \infty$.

Another quantity that seems to be more convenient if we want to 
compare two Markov chains is the spectral gap.
For this let the Markov chain $P$ be additionally reversible with 
respect to $\pi$, i.e.
\[
\pi(x)\,P(x,y) \;=\; \pi(y)\,P(y,x) \quad \text{ for all } x,y\in\O.
\]
(All the transition matrices from this article satisfiy this 
condition.)\\
Then we know that the eigenvalues of $P$ are real and we define the 
\emph{spectral gap} by
\[
\lambda(P)=1-\max\Bigl\{\abs{\xi}:\, \xi \text{ is an eigenvalue of } P,\; 
		\xi\neq1\Bigr\}.
\]
  
The eigenvalues of the Markov chain can be expressed in terms of norms 
of the operator $P$ that maps from 
$L_2(\pi):=(\R^\O,\pi)$ to $L_2(\pi)$, where 
inner product and norm are given by 
$\l f,g\r_\pi=\sum_{x\in\O}f(x) g(x) \pi(x)$ and 
$\Vert f\Vert_\pi^2:=\sum_{x\in\O}f(x)^2\pi(x)$, respectively. 
The operator is defined by
\begin{equation} \label{eq:map}
Pf(x) \;:=\; \sum_{y\in\O}\,P(x,y)\,f(y)
\end{equation}
and represents the expected value of the function $f$ after one step of 
the Markov chain starting in $x\in\O$. 
The \emph{operator norm} of $P$ is
\[
\Vert P\Vert_\pi \;:=\; \Vert P\Vert_{L_2(\pi)\to L_2(\pi)} 
\;=\; \max_{\Vert f\Vert_\pi\le1} \Vert Pf\Vert_\pi
\]
and we use $\Vert\cdot\Vert_\pi$ interchangeably for functions and 
operators, because it will be clear from the context which 
norm is used. It is well known that 
$\lambda(P)=1-\norm{P-S_\pi}_\pi$ for reversible $P$, 
where $S_\pi(x,y)=\pi(y)$. 
We know that reversible $P$ are self-adjoint with respect to 
the inner product 
$\l\cdot,\cdot\r_\pi$, i.e.~$P=P^*$, 
where $P^*$ is the \emph{adjoint operator} that satisfies 
$\l f, Pg\r_{\pi} = \l P^*f, g\r_{\pi}$ for all 
$f,g\in L_2(\pi)$.
The mixing time and spectral gap of a Markov chain (on finite state spaces) 
are closely related by the following inequality (see, 
e.g.,~\cite[Theorem~12.3 \& 12.4]{LPW}).

\begin{lemma} \label{lemma:mixing-gap}
Let $P$ be the transition matrix of a reversible, ergodic Markov chain with 
state space $\O$ and stationary distribution $\pi$. Then
\[
\lambda(P)^{-1}-1 \;\le\; \tau(P) 
\;\le\; \log\left(\frac{2e}{\pi_{\rm min}}\right)\,\lambda(P)^{-1}, 
\]
where $\pi_{\rm min}:=\min_{x\in\O}\pi(x)$.
\end{lemma}

\section{The algorithms} \label{sec:alg}

The \emph{Swendsen-Wang dynamics} (on the random-cluster model)
is based on the given connection of the random cluster and Potts models 
and performs the following two steps:
\begin{enumerate}
	\item[1)] Given a random cluster state $A\subset E$ 
		on $G$, assign a random color independently to each 
		connected component of $(V,A)$. Vertices of the same component 
		get the same color. This gives $\sigma\in\O_{\rm P}$.		
	\item[2)] Take $E(\sigma)$ and delete each edge independently with 
		probability $1-p$. This gives the new state $B\subset E$.
\end{enumerate}
This can be seen as first choosing $\sigma$ with respect to the 
conditional probability of $\bar\mu$ given $A$ and then 
choosing $B$ with respect to $\bar\mu$ given $\sigma$.
The transition matrix of the Swendsen-Wang dynamics is given by
\begin{equation} \label{eq:SW}
P_{\rm SW}(A,B) \;=\; q^{-c(A)}\,\left(\frac{p}{1-p}\right)^{\abs{B}}\,
	\sum_{\sigma\in\O_{\rm P}} \,(1-p)^{\abs{E(\sigma)}}\,
		\ind\bigl(\sigma\in\O(A\cup B)\bigr).
\end{equation}
Note that one can consider the Swendsen-Wang dynamics also on the Potts 
model (as it is usually done), which, starting at some $\sigma\in\O_{\rm P}$, 
performs the two steps above in reverse order. 
It is easy to prove that Swendsen-Wang on Potts and random-cluster model 
have the same spectral gap, see e.g.~\cite[Sec. 2.4]{U2}.

The second algorithm we want to analyze is the (lazy) \emph{single-bond dynamics}. 
Let $A\subset E$ be given and denote by $\stackrel{A}{\leftrightarrow}$ 
(resp.~$\stackrel{A}{\nleftrightarrow}$) connected (resp.~not connected) in 
the subgraph $(V,A)$.
Additionally we use throughout this article $A\cup e$ instead of $A\cup\{e\}$ 
(respectively for $\cap,\setminus$).
Note that $e_1\stackrel{A}{\nleftrightarrow}e_2$ for some $\{e_1, e_2\}\in E$ 
implies $\{e_1, e_2\}\notin A$.
The single-bond dynamics performs the following steps: 
\begin{enumerate}
	\item[1)] With probability $\frac12$ set $B=A$. 
		Otherwise, choose an edge $e=\{e_1,e_2\}\in E$ uniformly at random.		
	\item[2)] \begin{enumerate}
		\item[(i)] If $e_1\stackrel{A}{\leftrightarrow} e_2$:
			\begin{itemize}
				\item $B=A\cup e$ with probability $p$.
				\item $B=A\setminus e$ with probability $1-p$.
			\end{itemize}
		\item[(ii)] If $e_1\stackrel{A}{\nleftrightarrow} e_2$:
			\begin{itemize}
				\item $B=A\cup e$ with probability $\frac{p}{q}$.
				\item $B=A$\quad\;\; with probability $1-\frac{p}{q}$.
			\end{itemize}
		\end{enumerate}
	\item[3) ] The new state is $B$.
\end{enumerate}

The transition matrix of this Markov chain can be written as 
\begin{equation} \label{eq:SB}
P_{\rm SB}(A,B)
	\;=\; \frac{I(A,B)}{2}+\frac1{2\abs{E}}\sum_{e\in E}\,P_e(A,B),
\end{equation}
where $I(A,B)=\ind(A=B)$ and $P_e$ is given by
\begin{equation*} \label{eq:P-e}
P_e(A,B)
	\;=\; \ind(A\ominus B\subset e)\cdot\begin{cases}
	p^{\abs{B\cap e}}\,(1-p)^{1-\abs{B\cap e}}, 
		&e_1\stackrel{A}{\leftrightarrow} e_2\\
	\bigl(\frac{p}{q}\bigr)^{\abs{B\cap e}}\,(1-\frac{p}{q})^{1-\abs{B\cap e}}, 
		&e_1\stackrel{A}{\nleftrightarrow} e_2.
	\end{cases}
\end{equation*}
Here, $\ominus$ denotes the symmetric difference.

\begin{remark} \label{remark:HB}
If one is interested in the usual heat-bath dynamics on the random-cluster 
model, i.e. 
\begin{equation} \label{eq:P-HB}
\widetilde{P}(A,B) \;:=\; \frac1{2\abs{E}}\sum_{e\in E}\,
\frac{\mu(B)}{\mu(A\cup e)+\mu(A\setminus e)}\;\ind(A\ominus B\subset e) 
\quad\text{ for } A\neq B
\end{equation}
and $\widetilde{P}(A,A)$ such that $\widetilde{P}$ is stochastic, then
all results of this article hold up to a constant.
This is because 
$P_{\rm SB}(A,B)\le \widetilde{P}(A,B)\le (1-p(1-q^{-1}))^{-1}\,P_{\rm SB}(A,B)$ 
for all $A\neq B$ and so it is easy to prove by standard techniques 
(see e.g. \cite{DSC2}) that
\[
\lambda(P_{\rm SB}) \;\le\; \lambda(\widetilde{P}) 
	\;\le\; \bigl(1-p(1-q^{-1})\bigr)^{-1}\, \lambda(P_{\rm SB}).
\]
\end{remark}

\section{Representation on the joint model} \label{sec:joint}

In this section we want to represent the Swendsen-Wang and the single-bond 
dynamics on the FKES model, which consists of the product state space 
$\O_{\rm J}:=\O_{\rm P}\times\O_{\rm RC}$ and the FKES measure $\bar\mu$. 
For this we need the following ``building blocks''. 
First we introduce the stochastic matrix that defines the mapping (by matrix 
multiplication) from the RC to the FKES model
\begin{equation} \label{eq:M}
M\bigl(B,(\sigma,A)\bigr) \;:=\; q^{-c(B)}\;\ind\bigl(A=B\bigr)\;
	\ind\bigl(\sigma\in\O(B)\bigr).
\end{equation}
Note that $M$ defines an operator (like in \eqref{eq:map}) that maps from 
$L_2(\bar\mu)$ to $L_2(\mu)$ and its adjoint operator $M^*$ can be given 
by the stochastic matrix
\[
M^*\bigl((\sigma,A),B\bigr) \;=\; \ind\bigl(A=B\bigr).
\]
The following matrix represents the updates of the RC 
``coordinate'' in the FKES model. 
For $(\sigma,A),(\tau,B)\in\O_{\rm J}$ and $e=\{e_1,e_2\}\in E$ let
\begin{equation} \label{eq:T-e}
T_e\bigl((\sigma,A),(\tau,B)\bigr) \;:=\; \ind\bigl(\sigma=\tau\bigr)\;
	\begin{cases}
	p, & B=A\cup e  \,\text{ and }\;  \sigma(e_1)=\sigma(e_2)\\
	1-p, & B=A\setminus e  \;\text{ and }\;  \sigma(e_1)=\sigma(e_2)\\
	1, & B=A\quad  \;\,\;\text{ and }\;  \sigma(e_1)\ne\sigma(e_2)\\
	0, & \text{otherwise}.
	\end{cases}
\end{equation}
Clearly, one of the first two cases corresponds to $B=A$.
Since the transition probabilities $T_e\bigl((\sigma,A),(\tau,B)\bigr)$, 
which are 0 for $\sigma\neq\tau$, do not depend on whether $e\in A$ or not, 
it is convenient to state them in the above form.

Before we state the Swendsen-Wang and the single-bond dynamics in terms 
of the matrices from \eqref{eq:M} and \eqref{eq:T-e}, we state 
some properties that will be useful.

\begin{lemma} \label{lemma:prop}
Let $M$, $M^*$ and $T_e$ be the matrices from above.
Then
\begin{enumerate}
	\renewcommand{\labelenumi}{(\roman{enumi})}
	\item $M^*M$ and $T_e$ are self-adjoint in $L_2(\bar\mu)$.
\vspace{1mm}
	\item $M\,M^*(A,B)=\ind\bigl(A=B\bigr)$ and thus $M^*M\,M^*M = M^*M$.
\vspace{1mm}
	\item $T_e T_e = T_e$ and $T_e T_{e'} = T_{e'} T_e$ for all $e,e'\in E$.
\vspace{1mm}
	\item $\norm{T_{e}}_{\bar\mu}=1$ and $\norm{M^*M}_{\bar\mu}=1$.
\end{enumerate}
\end{lemma}
\begin{proof}
Part $(i)$ and $(ii)$ follow from the definition. 
Part $(iii)$ comes from the fact that the transition 
probabilities depend only on the ``coordinate'' that will not be changed and 
$(iv)$ follows from $(i)$, $(ii)$ and $(iii)$, since 
$\norm{T_{e}}_{\bar\mu} = \norm{T_{e}^2}_{\bar\mu}$ by $(iii)$ and 
$\norm{T_{e}^2}_{\bar\mu} = \norm{T_{e}}_{\bar\mu}^2$ by self-adjointness 
of $T_e$.\\
\end{proof}

Now we can state the desired Markov chains with the matrices from above. 

\begin{lemma} \label{lemma:repr}
Let $M$, $M^*$ and $T_e$ be the matrices from above.
Then
\begin{enumerate}
	\renewcommand{\labelenumi}{(\roman{enumi})}
	\item $P_{\rm SW} \,=\, M \left(\prod\limits_{e\in E} T_e\right) M^*$.
\vspace{1mm}
	\item $P_{\rm SB} \,=\, \frac{I}{2}
								+\frac1{2\abs{E}}\sum\limits_{e\in E}\,M\,T_e\,M^*$.
\end{enumerate}
\end{lemma}

From Lemma \ref{lemma:prop}$(iii)$ we have that the order of multiplication 
of the $T_e$'s 
in $(i)$ is unimportant.

\begin{proof}
For $(i)$ note that
\[
\biggl(\prod\limits_{e\in E} T_e\biggr)\bigl((\sigma,A),(\tau,B)\bigr) 
	\;=\; \ind\bigl(\sigma=\tau\bigr)\,\ind\bigl(B\subset E(\sigma)\bigr)\,
	p^{\abs{B}}(1-p)^{\abs{E(\sigma)}-\abs{B}}.
\]
Hence,
\[\begin{split}
M \left(\prod\limits_{e\in E} T_e\right) M^*\bigl(A,B\bigr) 
\;&=\; \sum_{\sigma\in\O_{\rm P}}\, M\bigl(A,(\sigma,B)\bigr)\,
	\biggl(\prod\limits_{e\in E} T_e\biggr)\bigl((\sigma,A),(\tau,B)\bigr)\\
&=\; \sum_{\sigma\in\O_{\rm P}} q^{-c(A)}\,
		\ind\bigl(\sigma\in\O(A)\cap\O(B)\bigr)\,
		p^{\abs{B}}(1-p)^{\abs{E(\sigma)}-\abs{B}}\\
&=\; P_{\rm SW}(A,B).
\end{split}\]
For $(ii)$ it is enough to prove $P_e=M\,T_e\,M^*$, where $P_e$ is from 
\eqref{eq:SB}. 
First we define 
$\ind_e(\sigma):=\ind\bigl(\sigma(e_1)=\sigma(e_2)\bigr)$ and 
$\ind_e(A):=\ind\bigl(e_1\stackrel{A}{\leftrightarrow} e_2\bigr)$ 
for $\sigma\in\O_{\rm P}$, $A\in\O_{\rm RC}$ and $e=\{e_1,e_2\}\in E$.
Now write 
\[
T_e\bigl((\sigma,A),(\sigma,B)\bigr) \;=\;
\ind\bigl(B=A\setminus e\bigr) + p\,\ind_e(\sigma)
\Bigl[\ind\bigl(B=A\cup e\bigr)
-\ind\bigl(B=A\setminus e\bigr)\Bigr]
\]
and note that $\abs{\O(A)}=q^{c(A)}$ and
\[
q^{-c(A)}\sum_{\sigma\in\O(A)}\ind_e(\sigma) 
\;=\; \frac1q \,+\, \ind_e(A)\left(1-\frac1q\right).
\]
Hence,
\[\begin{split}
M\,T_e\,M^*\bigl(A,B\bigr) 
\;&=\; \sum_{\sigma} q^{-c(A)}\ind\bigl(\sigma\in\O(A)\bigr)\,
				T_e\bigl((\sigma,A),(\sigma,B)\bigr)\\
&=\;\ind\bigl(B=A\setminus e\bigr) + p\,
		\Bigl[\ind\bigl(B=A\cup e\bigr)
		-\ind\bigl(B=A\setminus e\bigr)\Bigr]\cdot\\
		&\hspace{7cm}\cdot q^{-c(A)}\sum_{\sigma\in\O(A)}\ind_e(\sigma)\\
&=\; \begin{cases}
	p, & B=A\cup e  \,\text{ and }\; e_1\stackrel{A}{\leftrightarrow} e_2\\
	1-p, & B=A\setminus e \,\text{ and }\; e_1\stackrel{A}{\leftrightarrow} e_2\\
	\frac{p}{q}, 
		& B=A\cup e  \,\text{ and }\; e_1\stackrel{A}{\nleftrightarrow} e_2\\
	1-\frac{p}{q}, 
		& B=A\setminus e  \,\text{ and }\; e_1\stackrel{A}{\nleftrightarrow} e_2.
	\end{cases}\\
&=\; P_e(A,B) \quad 
	\text{ for all } A,B\in\O_{\rm RC} \text{ with } A\ominus B\subset e.
\end{split}\]
\end{proof}

\section{Main result} \label{sec:main}

In this section we prove the following theorem.

\begin{theorem} \label{th:main}
Let $P_{\rm SW}$ and $P_{\rm SB}$ be the transition matrices 
of the Swendsen-Wang and single-bond dynamics from \eqref{eq:SW} 
and \eqref{eq:SB}, respectively.
Then
\vspace{1mm}
\[
\lambda(P_{\rm SW}) \;\ge\; \lambda(P_{\rm SB}).
\vspace{1mm}
\]
This holds for arbitrary graphs $G$, $p\in(0,1)$ and $q\in\N$.
\end{theorem}

Before we prove the theorem we state some corollaries.
The first one gives an analogous inequality for the mixing times of the 
two algorithms. 

\begin{corollary} \label{coro:mix}
Let $P_{\rm SW}$ and $P_{\rm SB}$ be the transition matrices 
of the Swendsen-Wang and single-bond dynamics for the random-cluster model 
on $G=(V,E)$ with parameters $p\in(0,1)$ and $q\in\N$. Then
\[
\tau(P_{\rm SW}) \;\le\; 
	\left(3+\abs{E}\log\frac{1}{p(1-p)} + \abs{V}\log q\right)\,\tau(P_{\rm SB}).
\]
\end{corollary}
\begin{proof}
By Lemma \ref{lemma:mixing-gap} and Theorem \ref{th:main} we obtain
\[\begin{split}
\tau(P_{\rm SW}) \;&\le\; 
	\log\left(\frac{2e}{\mu_{\rm min}}\right)\,\lambda(P_{\rm SW})^{-1}
\;\le\; 
	\log\left(\frac{2e}{\mu_{\rm min}}\right)\,\lambda(P_{\rm SB})^{-1}\\
\;&\le\; 
	\log\left(\frac{2e}{\mu_{\rm min}}\right)\,\tau(P_{\rm SB}) + 1 \\
\;&\le\; \left(3+\log\left(\mu_{\rm min}^{-1}\right)\right)
	\,\tau(P_{\rm SB}).
\end{split}\]
Since $\mu_{\rm min}^{-1}$ can easily bounded by 
$\left(p(1-p)\right)^{-\abs{E}}\,q^{\abs{V}}$ the result follows.\\
\end{proof}

The next two corollaries show some rapid mixing results for the 
Swendsen-Wang dynamics. These are stated in terms of the spectral gap, but 
by Lemma \ref{lemma:mixing-gap} one can also use 
mixing times, adding the same factor as in Corollary \ref{coro:mix}. 
The first one is rapid mixing if the underlying graph is a tree, 
which is already known 
(see e.g. \cite{CF}, \cite{BCT}). The second is rapid 
mixing for graphs with bounded linear-width, which follows 
from a result of Ge and {\v{S}}tefankovi{\v{c}} \cite{GeS}.
For this we define the \emph{linear-width} of a graph $G=(V,E)$ as the 
smallest number $\ell$ such that there exists an ordering $e_1,\dots,e_{\abs{E}}$ 
of the edges with the property that for every $i\in[\abs{E}]$ there are at 
most $\ell$ vertices that have an adjacent edge in $\{e_1,\dots e_i\}$ and in 
$\{e_{i+1},\dots e_{\abs{E}}\}$. See \cite{GeS} for bounds on the 
linear-width of paths, cycles, trees, and in terms of a related quantity, 
the tree-width.

\begin{corollary} \label{coro:tree}
Let $P_{\rm SW}$ be the transition matrix 
of the Swendsen-Wang dynamics for the random-cluster model 
on a tree $T=(V,E)$. Then
\[
\lambda(P_{\rm SW})^{-1} \;\le\; 2\bigl(1-p(1-q^{-1})\bigr)^{-1}\,\abs{E}.
\]
\end{corollary}
\begin{proof}
Since $\mu_{p,q}^T$ is a product measure, we can write 
\[
P_{\rm SB}(A,B) \;=\; \frac{1}{\abs{E}}\sum_{e\in E} 
	\left(\frac{I+P_e}{2}\right)(A\cap e, B\cap e) \,\prod_{f\neq e} \ind(A\cap f = B\cap f),
\]
where $(I+P_e)/2$ can be seen (here) as a $2\!\times\!2$-matrix, 
i.e.~the transition matrix of the single-bond dynamics on a single edge. 
This matrix 
has the eigenvalues $1$ and $(1+p(1-q^{-1}))/2$.
We obtain by \cite[Lemma~12.11]{LPW} that 
$\lambda(P_{\rm SB})^{-1}=2\bigl(1-p(1-q^{-1})\bigr)^{-1}\,\abs{E}$.
This concludes the proof. \\
\end{proof}

Note that this bound improves the one given in \cite[Corollary~3.2]{BCT}, 
because it does not depend on maximum degree and depth of the tree.\\
The next results follows immediately from \cite{GeS}.

\begin{corollary} \label{coro:width}
Let $P_{\rm SW}$ and $P_{\rm SB}$ be the transition matrices 
of the Swendsen-Wang and single-bond dynamics for the random-cluster model 
on a graph $G=(V,E)$ with linear-width bounded by $\ell$. Then
\begin{equation} \label{eq:coro-width}
\lambda(P_{\rm SW})^{-1} \;\le\; \lambda(P_{\rm SB})^{-1} 
\;\le\; 4\abs{E}^2\,q^{\ell+1}.
\vspace{2mm}
\end{equation}
\end{corollary}
\begin{proof}
The first inequality is Theorem \ref{th:main} and the second follows from 
\cite{GeS}.
In this article the authors consider the Metropolis version of the 
single-bond dynamics. This Markov chain has transition probabilities 
\[
P_{\rm M}(A,A\ominus e) \;=\;	\frac1{2\abs{E}}\min\left\{1,\,q^{c(A\ominus e)-c(A)}
	\left(\frac{p}{1-p}\right)^{\abs{A\ominus e}-\abs{A}}\right\}, 
	\quad A\subset E,
\]
with $P_{\rm M}(A,A)$ such that $P_{\rm M}$ is a stochastic matrix. 
For this Markov chain they proof a lower bound on the congestion, 
which is defined as follows. Let $\G=\{\g_{AB}:\, A,B\subset E\}$, where 
$\g_{AB}$ are paths from $A$ to $B$ in the (directed) graph 
$\mathcal{H}=(\O_{\rm RC},\mathcal{E})$ with 
$\mathcal{E}=\left\{(A,B):\,P_{\rm M}(A,B)>0\right\}$. Then we define 
the \emph{congestion} of $P_{\rm M}$ (w.r.t. $\G$) by
\[
\varrho(P_{\rm M},\G) \;:=\; 
	\max_{(B_1,B_2)\in\mathcal{E}}\,\frac{1}{\mu(B_1)\, P_{\rm M}(B_1,B_2)}\,
	\sum_{A,C: (B_1,B_2)\in\g_{AC}} \abs{\g_{AC}} \mu(A)\, \mu(C),
\]
where $\abs{\g_{AC}}$ denotes the length of the path. 
The bound of \cite[Lemma~16]{GeS} is 
$
\varrho(P_{\rm M},\G) \le 2\abs{E}^2 q^\ell\; 
$
for a suitable choice of $\G$, and so we obtain by \cite[Prop.~1]{DS} 
(note that $P_{\rm M}$ is lazy) that
\[
\lambda(P_{\rm M})^{-1} \;\le\; 2\abs{E}^2\,q^\ell.
\]
It is easy to show that 
$P_{\rm M}(A,B)\le 2q \,P_{\rm SB}(A,B)$ for all $A,B\subset~\!\!\!E$. 
Thus, we can conclude by standard techniques 
(see e.g.~\cite[eq.~(2.3)]{DSC2}) that
\[
\lambda(P_{\rm SB})^{-1} \;\le\; 2q\, \lambda(P_{\rm M})^{-1} 
\;\le\; 4\abs{E}^2\,q^{\ell+1}.
\]
\end{proof}

\subsection{Proof of Theorem \ref{th:main}}

First we need the following technical lemma. 

\begin{lemma} \label{lemma:ineq}
Let $H$ and $G$ be two Hilbert spaces with corresponding 
inner products $\l f,f'\r_H$ and $\l g,g'\r_G$ for $f,f'\in H$ and 
$g,g'\in G$.
Furthermore, let $A:G\to H$ be a bounded linear operator with 
adjoint operator $A^*$, i.e. $\l A^*f,g\r_G=\l f,Ag\r_H$ for all 
$f\in H$, $g\in G$, and 
let $B:G\to G$ be a positive (i.e. $\l Bg,g\r_G\ge0$), bounded, 
self-adjoint linear operator. Then
\[
\norm{ABBA^*}_{H\to H} \;\le\; \norm{B}_{G\to G}\,\norm{ABA^*}_{H\to H},
\]
where $\norm{B}_{G\to G}:=\sup_{\norm{g}_G\le1} \norm{Bg}_G$ 
(resp.,~for $\norm{\cdot}_H$).
\end{lemma}

\begin{proof}
By the assumptions, $B$ has a unique positive square root $B^{\frac12}$, 
i.e. $B=B^{\frac12} B^{\frac12}$, which is again self-adjoint (see, 
e.g.,~\cite[Th. 9.4-2]{Krey}). We obtain
\[\begin{split}
\norm{ABBA^*}&_{H\to H} 
\;=\; \norm{AB}_{G\to H}^2 
\;\le\; \norm{A B^{\frac12}}_{G\to H}^2
	\norm{B^{\frac12}}_{G\to G}^2 \\
\;&=\; \norm{A B^{\frac12}\bigl(A B^{\frac12}\bigr)^*}_{H\to H}
	\norm{B^{\frac12}\bigl(B^{\frac12}\bigr)^*}_{G\to G} \\
\;&=\; \norm{ABA^*}_{H\to H} \norm{B}_{G\to G}.
\end{split}\]
\end{proof}

Now we are able to state the proof of Theorem \ref{th:main}.

\begin{proof}[Proof of Theorem \ref{th:main}]
Let $S_\mu$ and $S_{(\mu,\bar\mu)}$ be the operators that are induced 
(see~\eqref{eq:map}) 
by the (transition) matrices $S_\mu(A,B):=\mu(B)$ and 
$S_{(\mu,\bar\mu)}(B,(\sigma,A))=\bar\mu(\sigma,A)$ for all 
$A,B\subset E$,\, $(\sigma,A)\in\O_{\rm J}$. 
We get from Lemma \ref{lemma:prop}(iii) that 
$\prod_{e\in E} T_e = T_{e'}\prod_{e\in E} T_e$ for all $e'\in E$, 
and so
\begin{equation} \label{eq:Te1}
\prod_{e\in E} T_e \;=\; \prod_{e\in E} T_e \prod_{f\in E} T_f
\end{equation} 
and
\begin{equation} \label{eq:Te2}
\prod_{e\in E} T_e 
\;=\; \left(\frac{J}{2}+\frac{1}{2\abs{E}}\sum_{f\in E} T_{f}\right) 
			\prod_{e\in E} T_e,
\end{equation} 
where $J\bigl((\sigma,A),(\tau,B)\bigr):=\ind\bigl((\sigma,A)=(\tau,B)\bigr)$. 
Define the operators $\mathcal{T}:=\prod_{e\in E} T_e$, 
$T=\frac{1}{|E|}\sum_{e\in E}T_e$, and $N:=M-S_{(\mu,\bar\mu)}$.
It is easy to verify that 
$S_{(\mu,\bar\mu)} \mathcal{T} S^*_{(\mu,\bar\mu)} = S_{\mu}$ 
as well as 
$M \mathcal{T} S^*_{(\mu,\bar\mu)} 
= S_{(\mu,\bar\mu)} \mathcal{T} M^* = S_{\mu}$.
The same holds if $\mathcal{T}$ is replaced by $T$.
Thus, 
\begin{equation}\label{eq:N-M}
\begin{split}
\norm{N \mathcal{T} N^*}_\mu 
\;&=\; \norm{M \mathcal{T} M^* 
					- M \mathcal{T} S_{(\mu,\bar\mu)}^* 
					- S_{(\mu,\bar\mu)} \mathcal{T} M^*
					+ S_{(\mu,\bar\mu)} \mathcal{T} S_{(\mu,\bar\mu)}^*}\\
\;&=\; \norm{M \mathcal{T} M^* - S_\mu}_\mu,
\end{split}\end{equation}
and consequently 
\[\begin{split}
\norm{P_{\rm SW} - S_\mu}_\mu 
	\;&\stackrel{L.\text{\scriptsize \ref{lemma:repr}}}{=}\; 
			\norm{M \mathcal{T} M^* - S_\mu}_\mu 
	\;=\; \norm{N \mathcal{T} N^*}_\mu \\
	&\stackrel{\eqref{eq:Te1}}{=}\; 
			\norm{N \mathcal{T}	\mathcal{T} N^*}_\mu 
	\;=\; \norm{N \mathcal{T}}_{L_2(\bar\mu) \to L_2(\mu)}^2,
\end{split}\] 				
because $N \mathcal{T}$ 
induces an operator that maps from $L_2(\bar\mu)$ to $L_2(\mu)$ and the operator 
$\mathcal{T} N^*$ is 
its adjoint. Using submultiplicativity we obtain
\[\begin{split}		
\norm{P_{\rm SW} - S_\mu}_\mu \;&\stackrel{\eqref{eq:Te2}}{=}\; 
	\norm{N \left(\frac{J+T}{2}\right) 
				\mathcal{T}\,}_{L_2(\bar\mu) \to L_2(\mu)}^2 
	\;\le\; \norm{N \left(\frac{J+T}{2}\right)}_{L_2(\bar\mu) \to L_2(\mu)}^2 
				\norm{\mathcal{T}}_{\bar\mu}^2\\
	&\stackrel{L.\text{\scriptsize \ref{lemma:prop}}}{=}\; 
			\norm{N \left(\frac{J+T}{2}\right)}_{L_2(\bar\mu) \to L_2(\mu)}^2 
	\;=\; \norm{N \left(\frac{J+T}{2}\right)^2 N^*}_\mu \\
	&\stackrel{L.\text{\scriptsize \ref{lemma:ineq}}}{\le}\; 
			\norm{N \left(\frac{J+T}{2}\right) N^*}_\mu 
	\;=\;	\norm{M \left(\frac{J+T}{2}\right) M^*	- S_\mu}_\mu \\
	\;&\stackrel{L.\text{\scriptsize \ref{lemma:prop}}}{=}\; \norm{P_{\rm SB} - S_\mu}_\mu, 
\end{split}\] 
where the last inequality comes from Lemma \ref{lemma:ineq} with 
$H=L_2(\mu)$, $G=L_2(\bar\mu)$, 
$B=\frac{J+T}{2}$, and $A=N$. 
The next to last equality is proven in the same way as \eqref{eq:N-M}.
Note that 
$B$ is positive semidefinite (as a stochastic matrix with $B(x,x)\ge\frac12$ for all 
$x\in\O_{\rm J}$) and $\norm{B}_{\bar\mu}=1$.
This proves the claim, because 
$\lambda(P_{\rm SW})=1-\norm{P_{\rm SW} - S_\mu}_\mu $.\\
\end{proof}

\begin{remark}
Note that the proof of Theorem \ref{th:main} would be correct also in the 
case of the non-lazy single-bond dynamics, i.e.
$P_{\rm SB} \,=\,\frac1{\abs{E}}\sum_{e\in E}\,P_e$, 
because similar to \eqref{eq:Te2}, we have 
$\prod_{e\in E} T_e = \bigl(\frac{1}{\abs{E}}\sum_{f\in E} T_{f}\bigr) 
\prod_{e\in E} T_e$ and $\frac{1}{\abs{E}}\sum_{f\in E} T_{f}$ is 
positive semidefinite, since so are all $T_e$. 
But for convenience of the proof of Corollary 
\ref{coro:width} we choose to consider the lazy version.
\end{remark}

\section{Mixing time bounds on the torus} \label{sec:torus}

In this section we prove that the mixing time of the single-bond dynamics 
on the discrete $d$-dimensional torus of side length $L$ at the transition 
temperature is exponential in $L^{d-1}$, complementing a result of 
Borgs, Chayes and Tetali \cite{BCT}. 
For the upper bound we use the bound of Corollary \ref{coro:width} together 
with a bound of the linear-width of the $d$-dimensional torus. 
The lower bound follows from the lower bound of the mixing time of 
the Swendsen-Wang dynamics from \cite[Theorem~1.2]{BCT}.
Since this is closely related to (and also uses) their results, we refer to 
\cite{BCT} and the references cited therein for details. 
Let 
\[
T_{L,d} = (\Z/L\Z)^d
\]
be the $d$-dimensional torus of side length $L$. 
We will prove the following theorems.

\begin{theorem} \label{th:upper}
Let $P_{\rm SB}$ be the single-bond dynamics for the random-cluster model 
on $T_{L,d}$ with parameters $p\in(0,1)$ and $q\in\N$. Then for all $L,d\ge 2$ 
we have
\[
\tau(P_{\rm SB}) \;\le\; \exp\Bigl\{k_1(p) + k_2(q)\, L^{d-1}\Bigr\},
\]
where
\[
k_1(p) \;:=\; \log\left(1+\log\frac{1}{p(1-p)}\right)
\]
and
\[
k_2(q) \;:=\; 4 + 3\log q + \log(1+\log q).
\]
\end{theorem}

\begin{theorem} \label{th:lower}
Let $d\ge2$. Then there exists a constant $k_3=k_3(d)>0$ such that, 
for $q$ and $L$ sufficiently large, the single-bond dynamics for 
the random-cluster model on $T_{L,d}$ satisfies
\[
\tau(P_{\rm SB}) \;\ge\; \exp\Bigl\{k_3\,\beta_0\, L^{d-1}\Bigr\}
\quad \text{ for }\; p = 1-e^{\beta_0},
\vspace{1mm}
\]
where $\beta_0$ is the Potts transition temperature, i.e.
$\beta_0 = \frac1d \log{q} + \mathcal{O}(q^{-1/d})$.
\end{theorem}

First we prove Theorem \ref{th:lower}.

\begin{proof}
For the proof we have to consider the Swendsen-Wang dynamics for the 
Potts model with Gibbs measure $\pi=\pi_{\beta,q}^{T_{L,d}}$ 
(see section~\ref{sec:models}) for $\beta=\beta_0$, 
that performs the 
two steps of the Swendsen-Wang dynamics (as given in section~\ref{sec:alg}) 
in reverse order. We denote its transition matrix by $\widetilde{P}_{\rm SW}$.
These two algorithms have the same spectral gap if 
$p=1-e^{-\beta}$ (see \cite[sec.~2.4]{U2}). We obtain from Lemma 
\ref{lemma:mixing-gap} that
\[
\tau(P_{\rm SB}) +1 \;\ge\; \lambda(P_{\rm SB})^{-1}
\;\stackrel{Thm.\text{\scriptsize \ref{th:main}}}{\ge}\; \lambda(P_{\rm SW})^{-1}
\;=\; \lambda(\widetilde{P}_{\rm SW})^{-1} 
\;\ge\;\log\left(\frac{2e}{\pi_{\rm min}}\right)^{-1}\,\tau(\widetilde{P}_{\rm SW}),
\]
where $\pi_{\rm min}=\min_{\sigma\in\O_{\rm P}}\pi(\sigma)$.
Obviously, $\pi_{\rm min}\ge e^{-\beta\abs{E}} q^{-\abs{V}}$ for graphs $G=(V,E)$.
We know from Theorem 1.2 of \cite{BCT} that there exists a constant 
$k'_3=k'_3(d)>0$ such that, for $q$ and $L$ large enough, 
\[
\tau(\widetilde{P}_{\rm SW}) \;\ge\; 
\exp\Bigl\{k'_3\,\beta_0\, L^{d-1}\Bigr\} \quad \text{ for }\; \beta=\beta_0.
\]
Thus, since $\abs{V}=L^d$ and $\abs{E}=d L^d$ for $G=T_{L,d}$, we obtain
\[
\tau(P_{\rm SB}) +1 \;\ge\; \left(2+L^d\log q + \beta d L^d\right)^{-1} 
\exp\Bigl\{k'_3\,\beta_0\, L^{d-1}\Bigr\},
\]
which implies (again for $L$ large enough) that there exists a constant 
$k_3=k_3(d)>0$, such that
\[
\tau(P_{\rm SB}) \;\ge\; \exp\Bigl\{k_3\,\beta_0\, L^{d-1}\Bigr\}
\quad \text{ for }\; p = 1-e^{\beta_0}.
\]
\end{proof}

For the proof of Theorem \ref{th:upper} we need the following lemma.

\begin{lemma} \label{lemma:width}
The linear-width of $T_{L,d}$ does not exceed $\;2 L^{d-1}+1$.
\end{lemma}
\begin{proof}
For the proof we need three other ``widths'' of graphs, 
i.e.~\emph{path-width}, \emph{proper path-width}, and \emph{bandwidth}, 
but we omit their definition, because we do not need them here.
First note the following three facts:
\begin{enumerate}
	\item Linear-width is not larger than path-width+1 \cite[Lemma~2]{Fomin05}.
	\item Path-width is not larger than proper path-width \cite{Kaplan96}.
	\item Proper path-width equals bandwidth \cite[Theorem 3.2]{Kaplan96}.
\end{enumerate}
Therefore it is enough to prove that ${\rm bw}(T_{L,d})$; i.e.,~the bandwidth 
of $T_{L,d}$, is at most $2 L^{d-1}$. 
For this note that $T_{L,d}$ is the cartesian product of 
$d$ cycles $T_{L,1}$ of length $L$ and that ${\rm bw}(T_{L,1})=2$ 
(see \cite[Theorem~4.1.1]{Chinn82}). So we obtain by Corollary~4.3.2 of 
\cite{Chinn82} that ${\rm bw}(T_{L,d})\le2L^{d-1}$.
\end{proof}

Now we are able to prove Theorem \ref{th:upper}.

\begin{proof}[Proof of Theorem \ref{th:upper}]
Let $l$ be the linear-width of $T_{L,d}$. We know from 
Lemma~\ref{lemma:width} that $l\le2 L^{d-1}+1$, and so 
$l+1\le3 L^{d-1}$ since $L,d\ge2$. It follows from 
Corollary~\ref{coro:width} that 
\[
\lambda(P_{\rm SB})^{-1} \;\le\; 4 d^2 L^{2d} q^{3 L^{d-1}}.
\]
Set $\eta=\frac{1}{p(1-p)}$, and using Lemma~\ref{lemma:mixing-gap}, we obtain
\[\begin{split}
\tau(P_{\rm SB}) 
\;&\le\; \log\left(\frac{2e}{\mu_{\rm min}}\right) \lambda(P_{\rm SB})^{-1}\\
&\le\; \left(2+L^d\log q + d L^d\log\eta\right) 4 d^2 L^{2d} q^{3 L^{d-1}} \\
&\le\; 4 d^3 L^{3d} q^{3 L^{d-1}} \left(1+\log q + \log\eta\right) \\
&=\; \exp\Bigl\{\log(4 d^3 L^{3d}) + 3\log(q) L^{d-1} 
			+ \log(1 + \log q + \log\eta)\Bigr\} \\
&\le\; \exp\Bigl\{4 L^{d-1} + 3\log(q) L^{d-1} 
			+ \log(1 + \log q) + \log(1 + \log\eta)\Bigr\} \\
&\le\; \exp\Bigl\{k_1(p) + k_2(q) L^{d-1}\Bigr\}
\end{split}\]
with $k_1$ and $k_2$ from Theorem~\ref{th:upper}. This proves the claim.
\end{proof}

{

\end{document}